\documentclass[12pt]{article}
\usepackage{graphicx}
\usepackage{amssymb}
\usepackage{amsmath}
\usepackage{cite}
\numberwithin{equation}{section}
\def\BC{{\mathbb C}}
\def\BR{{\mathbb R}}
\def\clp{{\mathcal P}}

\def\erf{{\rm erf\ }}
\def\col{{\rm col\ }}

\def\tr{{\rm tr\ }}
\newtheorem{Pa}{Paper}[section]
\newtheorem{Tm}[Pa]{{\bf Theorem}}
\newtheorem{La}[Pa]{{\bf Lemma}}
\newtheorem{Cy}[Pa]{{\bf Corollary}}
\newtheorem{Rk}[Pa]{{\bf Remark}}

\newtheorem{Ee}[Pa]{{\bf Example}}
\newtheorem{Dn}[Pa]{{\bf Definition}}
\newtheorem{Pn}[Pa]{{\bf Proposition}}

\def\Xint#1{\mathchoice
   {\XXint\displaystyle\textstyle{#1}}%
   {\XXint\textstyle\scriptstyle{#1}}%
   {\XXint\scriptstyle\scriptscriptstyle{#1}}%
   {\XXint\scriptscriptstyle\scriptscriptstyle{#1}}%
   \!\int}
\def\XXint#1#2#3{{\setbox0=\hbox{$#1{#2#3}{\int}$}
     \vcenter{\hbox{$#2#3$}}\kern-.5\wd0}}

\def\dashint{\Xint-}

\title{Characteristic function of M. Liv\v{s}ic \\ and  some developments}

\author{Lev Sakhnovich}

\date{}

\begin{document}

\maketitle

{\centering \it
To the memory of my teacher M. Liv\v{s}ic with admiration 
and gratitude. \par}

\begin{abstract}
The area  related  to  M. Liv\v{s}ic's characteristic matrix functions is too vast to be discussed in one paper
and we selected for this article the problems which are close to our scientific interests.
We  discuss M.Liv\v{s}ic's results connected with characteristic matrix functions and various important developments including
factorization of the transfer operator function, inverse problems, triangular models, reduction of the operators to the simplest
form as well as applications to Wiener--Masani prediction theory, to random matrices, and to Riemann--Hilbert problems.

\end{abstract}

\section{Introduction}
In this paper, we  focus on the famous article  \cite{Liv2} by M. Liv\v{s}ic, where two fundamental problems have been solved:

1. The notion of the characteristic  matrix-valued function (matrix function) $W_A(z)$ of the linear operator $A$
(acting in the  Hilbert space $H$) was introduced. It was proved that this
characteristic  matrix function  $W_A(z)$ determines the operator $A$ up to unitary equivalence.

2. Under the assumption that
the operator $(A-A^{*})/i$ is compact and its eigenvalues $z_n$ satisfy the inequality
\begin{equation}\sum_{n=1}^{\infty}|z_n|<\infty,\label{1.1}\end {equation}
the triangular model of the operator $A$ was constructed.
(M. Liv\v{s}ic  considered this result as the main one in his article, see \cite{Liv2}.)

In his works, M. Liv\v{s}ic studied operators $A$ close to self-adjoint operators and corresponding characteristic  matrix functions $W_A(z)$.
Often, we omit index
$A$  and write $W(z)$. We use also the notation $W_A(z)$ for the transfer matrix functions which generalize characteristic  matrix functions
(see section~3). Matrices and operators
$$J=J^*=J^{-1}$$
and the so called $J$-properties of $W(z)$ play an essential role in the theory.

\begin{Rk}\label{Remark 1.1}The famous Liv\v{s}ic's work \cite{Liv2} is closely connected with the famous Potapov's paper \cite{Pot}. These papers were written almost simultaneously with the mutual influence and cooperation of Liv\v{s}ic and Potapov. I am proud that I could read both of these articles $($when I was a student$)$ even before they were published.\end{Rk}
Here,  we describe some results where  M. Liv\v{s}ic's  ideas and assertions  have been used and developed.
Since the area is huge, we selected the results which are close to our scientific interests.

In section 2, we introduce the notion of the characteristic matrix function $W(z)$
and formulate the fundamental result on the unitary  equivalence of operators (M. Liv\v{s}ic \cite{Liv2}). It was further shown that the characteristic matrix function $W(z)$
plays an important role in the problems of operator similarity \cite{NF, Sak7, Sak11}.
M. Liv\v{s}ic commented with pleasure on the corresponding results  \cite{NF, Sak7, Sak11}:\\
``Thus, the characteristic matrix function does not only define the unitary invariants of an operator, but is also useful  in the problems of  linear similarity."

In section 3, we consider   factorization of the operator-valued transfer function $W(z)$. We note that the  characteristic matrix function is a partial case of  the operator transfer function. The method of factorization of the characteristic matrix function is given in the Liv\v{s}ic--Potapov paper \cite{LiPo}.
We understood  long ago that many formulas of the characteristic matrix function theory are purely algebraic relations and that they can be generalized for matrix functions
 without $J$-properties. It was not easy and took a long time, but we finally found the corresponding approach.  A method,  which permits  to analyse in one uniform manner various factorization problems, was proposed in \cite{Sak13, Sak2}. This method is described in section 3.

Under certain conditions, the factorization theorem (section 3) permits us
to find  a continual factorization of an operator function $W(z)$, that is, to
represent $W(z)$ as a multiplicative integral
\begin{equation}W(z)=\int\limits_a^{\overset{b}{\curvearrowleft}} \, \exp\left[ -\frac{d\sigma(t)}{\alpha(t)-z}\right]. \label{1.2}
\end{equation} V.P. Potapov  proved \cite{Pot} representation \eqref{1.2} for a much broader class. In \cite{Pot}, the existence problem for representation \eqref{1.2} was fully and completely solved.  In section 4, we describe   the method of construction of $\sigma(t)$ given in \cite{Sak2}. For this purpose,
we consider the limits
\begin{equation}W_{\pm}(x)=\lim{W(z)}\quad (z=x+iy, \,\, y{\to}\pm{0}). \label{1.3}\end{equation}
Then, we write down the polar decomposition
\begin{equation}W_{+}(x)=U(x)R(x),     \label{1.4}\end{equation}
where the operator functions $U(x)$ and $R(x)$ are such that
\begin{equation}U^{*}(x)JU(x)=J,\quad JR(x)=R^{*}(x)J, \label{1.5-}\end{equation}
and the spectrum of $R(x)$ is positive. The operator function $R(x)$ is the $J$-module of $W(x)$ \cite{Pot}.
 It is essential that the polar decomposition for $W_-(x)$ has the form \cite{Sak7, Sak2}
\begin{equation}W_{-}(x)=U(x)R^{-1}(x).     \label{n4}\end{equation}
In section 4, we solve the problem of recovering $W(z)$ from $J$
and from the given $J$-module $R(x)$.

In section 5, we  consider separately the special case where $J=I$. This case plays an important role in the Masani--Wiener prediction theory \cite{WiMa}. We consider  the problem to recover $W(z)$ from the given $I$-module $R(x)$, and prove  that the constructed  $W(z)$ is  the unique maximal solution.

In section 6 we investigate a special class of $J$-modules $R(x)$. Namely,
 we assume that $R(x)$ satisfies the following conditions\\
 \begin{equation}[R(x)-I]^{2}=0,\quad  m=2 ,\label{1.5}\end{equation}
 and
 \begin{equation}J=\left(
                     \begin{array}{cc}
                       -1 & 0 \\
                       0& 1 \\
                     \end{array}
                   \right).\label{1.6}\end{equation}
\begin{Dn}\label{Definition 1.2} By Z-class (zero-class) we denote the set of $J$-modules $R(x)$ which satisfy the condition \eqref{1.5}.\end{Dn}
In section 6, we use  Z-class theory  for the study of
the following important in random matrix theory cases: sine kernel (see \cite {DIZ}), Bessel kernel (see \cite{TrWi}) and Airy kernel (see \cite{TrWi1}).
In section 7, we study the connection between Fredholm determinant
and Riemann--Hilbert problem. In the sine-kernel case, we obtain a well--known
result from
\cite{DIZ}. Section 8 is dedicated to  the reduction of operators to
the triangular form.
The first important result in this area was obtained by M. Liv\v{s}ic,
and  he considered it as the main achievement in
\cite{Liv2}. In section 8, we formulate several assertions from \cite{Sak15, Sak16}
which generalize further the  results by  M. Livsic.
These new results served as an  impetus for a number of studies including research on the
existence of invariant subspaces for operators of the form  $A=H+iB$,
where B is a compact operator
\cite{GoKr1, GoKr2, GoKr4, Mac, Sak16},
and research on the interconnections between the imaginary and real parts of a  compact operator
with spectrum concentrated at zero (see \cite{GoKr1, GoKr2, GoKr4, Mac, Sak16} as well).

In section 9, we consider a simple but non-trivial example of a triangular
Liv\v{s}ic model and reduce it to the diagonal form \cite{Sak10}. The reduction
of operators to the simplest form  is a vast topic, and it cannot
be described it in this comparatively short paper.
Section 10 is auxiliary and is devoted to Christoffell-Darboux type formulas

{\it Notations.}
As usual $\BR$ stands for the real axis, $\BC$ stands for the complex plane, $\BC_+$ stands for the open upper half-plane,
and $\overline{z}$, where $z\in \BC$, denotes the complex conjugate
of $z$.
The notation $\overline{G}$ denotes the closure
of the space $G$, and
the symbol $[H_1,H_2]$
denotes the set of bounded operators acting from $H_1$ into $H_2$. The space of square integrable column vector functions
$f(x)\in \BC^N$ on $(a,b)$ is denoted by $L^2_N(a,b)$ and $L^2_1(a,b)=L^2(a,b)$.
The space of integrable on $(a,b)$ functions is denoted by $L(a,b)$.
The trace of some operator $K$ is denoted
by $\tr(K)$ and $K^*$ stands for the adjoint operator. The notations $\Re$ and $\Im$ of the real and imaginary parts
is extended to the bounded operators in Hilbert spaces: $\Re(A)=(A+A^*)/2$ and $\Im(A)=(A-A^*)/(2i)$.
\section{Characteristic matrix function}
The notion  of the
characteristic matrix function $W(z)$ (of the operator $T$ acting
in a separable Hilbert space $H$) was initially introduced in \cite{Liv2}. This notion was
further developed (and somewhat changed) by M. Liv\v{s}ic himself and by many other authors (see \cite{Liv2, BL}).
We shall use the following definition of the characteristic operator  function (of the operator $T$ which is bounded together with its inverse in $H$). 
\begin{Dn}\label{Definition 2.1} \cite{Sak7}.
Let a bounded operator $R$  map Hilbert space $G$ into Hilbert space $H$  and satisfy the conditions\\
1. $I-T^{*}T=RJR^{*}$, where $J$ acts in $G$ and $J=J^{*},\,\, J^2=I.$\\
2. The spectrum of the operator $I-JRR^{*}$ is strictly positive.\\
 Condition 2 implies that  there is a $J$-Hermitian operator $W(0)$ which has a strictly positive
spectrum and  the following property:
\begin{equation}W^{2}(0)=I-JRR^{*}.\nonumber\end{equation}
Now, the characteristic operator function $W(\mu)$ of the operator T is introduced by the
equality:
\begin{equation}W(0)W(\mu)=I-JR^{*}(I-\mu{T})^{-1}R.\label{2.1}\end{equation}\end{Dn}
M. Liv\v{s}ic obtained the following classical result \cite{Liv2}.
\begin{Tm}\label{Theorem 2.2} If two operators have the same characteristic matrix function,
then the simple parts of these operators are unitary equivalent. \end{Tm}
\begin{Rk}\label{Remark 2.3} Recall that the simple part of $T$ is the operator induced  by $T$ on the subspace\\
$$H_1=\overline{\sum_{k=-\infty}^{\infty}T^{k}D_T},\quad where \quad D_T=\overline{(I-T^{*}T)H}.$$
\end{Rk}
B. Sz.-Nagy and C. Foias \cite{NF}  proved the following assertion.
\begin{Tm}\label{Theorem 2.4} If $J=I$ then the boundedness condition
\begin{equation} \sup \|{W(\mu)}\|{\leq}C<\infty \quad for \quad |\mu|{\ne}1 \label{2.2}\end{equation}
is necessary and sufficient  for the simple part of $T$ to be linearly similar to a unitary operator with absolutely continuous spectrum.\end{Tm}
Later we  proved \cite{Sak7, Sak11}
 the following result  without assumption $J=I$.
\begin{Tm}\label{Theorem 2.5} Let the characteristic operator function $W(\mu)$
of the operator $T$ satisfy the condition \eqref{2.2}.
Then, the simple part of $T$ is linearly similar to a unitary operator with absolutely continuous spectrum.\end{Tm}
\begin{Rk}\label{Remark 2.6}  Theorems \ref{Theorem 2.4} and  \ref{Theorem 2.5} deal with the operators linearly similar to unitary operators. These theorems are easily   reformulated
for the case of operators linearly similar to self-adjoint operators (see \cite{Sak14}).\end{Rk}
\begin{Rk}\label{Remark 2.7} In \cite{Sak14}, Theorem \ref{Theorem 2.5} was applied to the analysis of partial differential operators. \end{Rk}
The next section is dedicated to the factorization theorem for  the operator-valued transfer function $W(z)$, which generalizes
corresponding results for the characteristic operator function.
Under certain conditions, this factorization theorem (section 3) permits us to find a continual
factorization of $W(z)$, that is, to represent $W(z)$ as a multiplicative integral (see section 4).

\section{On the factorization of the operator-valued transfer function}
Let us consider a linear dynamical system
\begin{align} & dx/dt=Ax+\Pi_{1}v, \quad y=\Gamma_{2}^{*}x+v;\label{3.1}
\\ & \label{n1}
\Pi_{1}{\in}[G,H_1],\quad A{\in}[H_1,H_1],\quad \Gamma_{2}^{*}{\in}[H_1,G],
\quad  x(t){\in}H_1,\quad  y(t),\, v(t){\in}G.
\end{align}
where $H_1$ and $G$ are  Hilbert spaces. It is well known that the operator .transfer function $W_A(z)$ of system \eqref{3.1} has the form
\begin{equation}W_A(z)=I-\Gamma_{2}^{*}(A-zI)^{-1}\Pi_{1}.\label{3.2}\end{equation}
A representation of a given operator function in the form \eqref{3.2} is called  {\it realization}.We  use also a realization of $W_{A}^{-1}(z)$:
\begin{equation}W_{A}^{-1}(z)=W_B(z)=I+\Pi_{2}^{*}(B-zI)^{-1}\Gamma_{1},\label{3.3}\end{equation}
where
\begin{align} &\label{n2}
\Gamma_{1}{\in}[G,H_2],\quad B{\in}[H_2,H_2],\quad\Pi_{2}^{*}{\in}[H_2,G],
\end{align}
and $H_2$ is a Hilbert space.
Without loss of generality, we may assume \cite{Kal} that the realizations \eqref{3.2} and \eqref{3.3} are minimal, that is,
\begin{equation}\overline{\sum_{n=0}^{\infty}A^{n}\Pi_{1}{G}}=H_{1},\quad
\overline{\sum_{n=0}^{\infty}(A^*)^{n}\Gamma_{2}{G}}=H_{1},\label{3.4}\end{equation}
\begin{equation}\overline{\sum_{n=0}^{\infty}B^{n}\Gamma_{1}{G}}=H_{2},\quad
\overline{\sum_{n=0}^{\infty}(B^{*})^{n}\Pi_{2}{G}}=H_{2}.\label{3.5}\end{equation}
We introduce also an important operator $S{\in}[H_1,H_2]$ satisfying the equation
\begin{equation}        AS-SB=\Pi_{1}\Pi_{2}^{*},              \label{3.6}\end{equation}
and the conditions
\begin{equation}S\Gamma_1=\Pi_1,\quad \Gamma_{2}^{*}S=\Pi_{2}^{*}.\label{3.7}\end{equation}
It follows from \eqref{3.6} and \eqref{3.7} that
\begin{align}   &     TA-BT=\Gamma_{1}\Gamma_{2}^{*},             \label{3.8}
\\ &
T\Pi_1=\Gamma_1,\quad \Pi_{2}^{*}T=\Gamma_{2}^{*}S,\label{3.9}\end{align}
where
\begin{equation}T=S^{-1}.\label{3.10}\end{equation}
The operator identities \eqref{3.6}-\eqref{3.10} are basic in the method of operator identities and play an essential role in our considerations.
Further in the text, we assume that operators $S$ and $T$ are bounded and
\begin{equation}H_1=H_2=H.\label{3.11}\end{equation}
Now, let the space $H$ be represented in the form
\begin{equation}H=L_1{\oplus}L_2.     \label{3.12}\end{equation}
The operator of orthogonal projection of $H$ onto $L_k$ is  denoted by $P_k$.
We proved the following fundamental factorization theorem \cite{Sak13, Sak2}.
\begin{Tm}\label{Theorem 3.1} Let the operator-valued transfer function $W_{A}(z)$
satisfy  conditions \eqref{3.2}--\eqref{3.5},
let  the bounded operator S satisfy  \eqref{3.6} and \eqref{3.7}, and let  \eqref{3.11} and \eqref{3.12} hold.
Assume that the operators A and B have  triangular forms, that is,
\begin{equation}AP_2=P_{2}AP_{2},\quad BP_1=P_{1}BP_{1}.\label{3.13}\end{equation}
If the operators S and $S_{11}=P_{1}SP_{1}$ are invertible on H and $L_1$, respectively,
then the following formula is valid:
\begin{equation}W_A(z)=W_2(z)W_1(z), \label{3.14}\end{equation}
where
\begin{equation}W_1(z)=I-\Pi_{2}^{*}P_{1}S_{11}^{-1}(A_{11}-zI)^{-1}P_{1}\Pi_{1},\label{3.15}\end{equation}
\begin{equation}W_2(z)=I-\Gamma_{2}^{*}P_{2}(A_{22}-zI)^{-1}T_{22}^{-1}P_{2}\Gamma_{1},\label{3.16 }\end{equation}
\begin{equation}A_{jj}=P_{j}AP_{j},\quad T_{22}=P_2S^{-1}P_2.\label{3.17}\end{equation}\end{Tm}
\begin{Rk}\label{Remark 3.2}The invertibility of $T_{22}$ on the subspace $L_{2}$ stems  from the invertibility of S and $S_{11}$ on H and $L_1$, respectively.\end{Rk}
Theorem \ref{Theorem 3.1} is close in spirit to the earlier Liv\v{s}ic--Potapov theorem   \cite{LiPo} on the factorization of the characteristic matrix function. In order
to explain the connection between these theorems we write $A-A^{*}$ in the form
\begin{equation} A-A^{*}=i\Pi_{1}J\Pi_{1}^{*},\quad J=J^{*},\quad J^{2}=I,\quad J{\in}[G,G].\label{3.18}\end{equation}
Relation \eqref{3.18} is a special case of \eqref{3.6} where
\begin{equation} S=I,\quad B=A^{*},\quad \Pi_{2}^{*}=iJ\Pi_{1}^{*}.\label{3.19}\end{equation}
In this case,  the transfer operator function $W_A(z)$ defined by \eqref{3.2} has the form
\begin{equation}W_A(z)==I-J\Pi_{1}^{*}(A-zI)^{-1}\Pi_{1}\label{3.20}\end{equation}
and coincides with the characteristic operator function of the operator A \cite{Liv2}. Thus,
the well-known theorem on factorization of the characteristic operator function \cite{LiPo}
follows from Theorem \ref{Theorem 3.1}.

A geometrical interpretation and applications of Theorem 3.1 are contained in the book
\cite{BGK}. A version of Theorem 3.1 related to quasi-unitary operators was obtained in
\cite{Zol}. Theorem 3.1 found also a  number of applications in our works (see section 4).
\section{Continual factorization \\ and inverse problems}
Under certain conditions the factorization theorem (section 3) permits us to find a continual
factorization of an operator function $W(z)$, that is to represent $W(z)$ as a multiplicative integral
\begin{equation}
W(z)=\int \limits_a^{\overset{b}{\curvearrowleft}} \, \exp\left[ -\frac{d\sigma(t)}{\alpha(t)-z}\right]. \label{4.1}
\end{equation}
The construction of $\sigma(t)$ presented in our paper \cite{Sak2}  yields a unified solution
for a wide range of  problems in several domains:\\
1. The construction of Liv\v{s}ic triangular model which corresponds to characteristic
operator function $W(z)$.\\
2. Prediction problems for stationary processes (see Wiener and Masani \cite{WiMa}, the case $\alpha(x)=x)$.\\
3. Riemann--Hilbert problem.\\
4. Random matrix theory. \\
5. Classical inverse problems, for instance, for radial Schr\"odinger equations and  Dirac type systems
(the case $\alpha(x)=0 $, see \cite{Sak8}).\\
6. Integrable operators \cite{Sak2}.

Let us consider in detail the case where $\alpha(x)=x$, and the matrix function $W(z)$ is acting in $G$ for each $z$ and is analytic
in the domain $z{\ne}\overline{z}$. We suppose that $\dim\,{G}=N<\infty$ and set (without loss of generality) $G=\BC^N$.
We assume  the following conditions are valid for $W$.\\
1. The $N\times N$ matrix function  $W(z)$  satisfies the condition
\begin{equation}\sup \|W(z)\|{\leq}M,\quad z{\ne}\overline{z}.\label{4.2}\end{equation}
2. All the points $z\in\BC$ excluding some finite interval $[a,b]\subset \BR$ $(z{\notin[a,b]})$ are regular points of $W(z)$ and
\begin{equation} \lim W(z)=I\quad (z{\to}\infty) .\label{4.3}\end{equation}
3. The inequality
\begin{equation}i\frac{W(z)JW^{*}(z)-J}{z-\overline{z}}{\geq}0
\label{4.4}\end{equation}
holds for $z{\ne}\overline{z}$ and for $J=J^{*}=J^{-1}$.
\begin{Rk}\label{Remark 4.1} Each characteristic matrix function $W(z)$ satisfies conditions
\eqref{4.3} and \eqref{4.4}.  Condition \eqref{4.2} means that the corresponding operator A
is linearly similar to some self-adjoint operator with absolutely continuous spectrum \cite{Sak14}.\end{Rk}

Now, let us construct operators  $\Pi_1$, $\Gamma_1$, $\Gamma_2^*$ and $\Pi_2^*$ acting
in the spaces given in \eqref{n1}, \eqref{n2}, where $H_1=H_2=H=L^2_N(a,b)$, so that the realizations \eqref{3.2} and \eqref{3.3} hold.

Taking into account conditions   \eqref{4.2}--\eqref{4.4} and using Fatou theorem, one can show that the matrix functions  $R(x)$ and $U(x)$,
which  give the polar decomposition \eqref{1.4} and satisfy  \eqref{1.5-} are unique and exist almost everywhere.
Moreover, we have the following lemma.
\begin{La}\label{Lemma 4.3}\cite{Sak2}. Under conditions \eqref{4.2}--\eqref{4.4}, the matrix-valued functions $R(x)$, $R^{-1}(x)$,
$U(x)$, and $U^{-1}(x)$ are bounded on [a,b] and the next inequality holds:
\begin{equation}D(x)=J[R(x)-R^{-1}(x)]{\geq}0.\label{4.10}\end{equation}\end{La}
Introduce  matrix functions $F_{1}(x)$, $F_{2}(x)$, and $L(x)$ such that
\begin{align} & F_{1}^{*}(x)F_{1}(x)=D(x),\quad F_{2}(x)=F_{1}(x)JU^{*}(x),\label{4.11}
\\ &
L(x)=\{I+\frac{1}{4}[F_{1}^{*}(x)JF_{1}(x)]^{2}]\}^{1/2}.
\label{4.12}\end{align}
Now, the operators $\Pi_{1}$ and  $\Gamma_{2}^*$ are given by the equalities
\begin{equation} \Pi_{1}g=\frac{1}{\sqrt{2\pi}}F_{1}(x)g,\quad  \Gamma_{2}^{*}f(x)=\frac{i}{\sqrt{2\pi}}\int_{a}^{b}F_{2}^{*}(x)f(x)dx.
\label{4.14}\end{equation}
Note that
\begin{equation} \Pi_{1}^{*}f(x)=\frac{1}{\sqrt{2\pi}}\int_{a}^{b}F_{1}(x)f(x)dx,
\quad \Gamma_{2}g=-i\frac{1}{\sqrt{2\pi}}F_{2}(x)g,\quad g{\in}G.
\label{4.15}\end{equation}
The operators $\Pi_2$ and $\Gamma_{1}$ are expressed via $\Pi_{1}$ and  $\Gamma_{2}$:
\begin{equation}\Pi_{2}=-i\Pi_{1}J,\quad \Gamma_{1}=i\Gamma_{2}J, \label{4.18}\end{equation}
and the operators $A$ and $B$ have the form
\begin{equation} Af=Bf=xf(x),\quad f(x) \in  L^2_N(a,b).\label{4.19}\end{equation}
Then, the realizations
\begin{align} & W(z)=I-\Gamma_{2}^{*}(A-zI)^{-1}\Pi_{1},\label{4.16}
\\ &
W^{-1}(z)=I+\Pi_{2}^{*}(A-zI)^{-1}\Gamma_{1},\label{4.17}\end{align}
are valid.

The orthogonal projector from $L^2_N(a,b)$ onto $L^2_N(a,\zeta)$ $(a<\zeta \leq b)$ is denoted
by $P_{\zeta}$.
We write down  the operators
\begin{align} &
S_{\zeta}=P_{\zeta}S_b P_{\zeta} \, \in \, [L^2_N(a,\zeta), L^2_N(a,\zeta)]  \quad(a<\zeta \leq b)
\label{4.17+}\end{align}
and $S_{\zeta}^{-1}$ in the following way:
\begin{align} & S_{\zeta}f=L(x)f(x)+\frac{i}{2\pi}\dashint_a^\zeta\frac{F_{1}^{*}(x)JF_{1}(t)}{x-t}f(t)dt,
\label{4.20}
\\  &
S_{\zeta}^{-1}f=L(x)f(x)-\frac{i}{2\pi}\dashint_a^\zeta\frac{F_{2}^{*}(x)JF_{2}(t)}{x-t}f(t)dt,
\label{4.21} \end{align}
where $f\in L^2_N(a,\zeta)$ and $\dashint_a^b$ denotes Cauchy integral.

Now, we are ready to formulate a fundamental multiplicative representation theorem \cite[Ch.III, section 2]{Sak2}.
\begin{Tm}\label{Theorem 4.2} Let conditions \eqref{4.2}--\eqref{4.4} be fulfilled.
Then, $W(z)$ admits multiplicative representation
\begin{equation}
W(z)=\int \limits_a^{\overset{b}{\curvearrowleft}} \, \exp\left[ -\frac{d\sigma(t)}{t-z}\right],\label{4.5}
\end{equation}
where
\begin{equation}\sigma(\zeta)=\Pi_{2}^{*}S_{\zeta}^{-1}P_{\zeta}\Pi_{1},\label{4.6}\end{equation}
$\Pi_1$ is given in \eqref{4.14}, $\Pi_2$ is given in \eqref{4.18}, and $S_{\zeta}^{-1}$ is given in \eqref{4.21}.
\end{Tm}
\begin{Rk}\label{Rk4.4} According to Theorem \ref{Theorem 4.2}, it suffices to know  $R(x)$ and $J$ in order to construct
$\sigma(\zeta)$.   \end{Rk}
The problem of recovering $W(z)$ from the given
$R(x)$ and $J$ is of independent interest and has applications to prediction problems for stationary processes (the case case $J=I$, see \cite{Ros,WiMa})) and to
Riemann--Hilbert problems. Note that we not only prove  the existence of  $W(z)$ but describe above two ways to recover $W(z)$. That is, $W(z)$ may be recovered either using realization \eqref{4.16} or multiplicative
representation \eqref{4.5}, \eqref{4.6}.

Application of the above-mentioned results to  prediction theory is discussed in the next section. Recall also  the Riemann-Hilbert problem:
\emph{if $R(x)$ and $J$ are given,  find $W(z)$ such that
\begin{equation}W_{+}(x)=W_{-}(x)R^{2}(x),\quad a{\leq}x{\leq}b.\label{4.22}\end{equation}}
It follows  from  \eqref{1.3}, \eqref{1.4}, and \eqref{n4} that
$W(z)$ recovered from $R(x)$ using \eqref{4.16} or \eqref{4.5}, \eqref{4.6} gives a solution of the problem \eqref{4.22}.
\begin{Rk}\label{Remark 4.5} If the elements of R(x) are rational, then the method of construction of  $W(z)$ is well known (see, e.g., \cite{Ros}).\end{Rk}
\begin{Rk}\label{4.9} The limiting values  $W_{\pm}(x)$ of a multiplicative integral $W(z)$ are found in the papers \cite{Sak4, Sak1, Sak6}.
The corresponding result can be regarded as an extension (to  the multiplicative integral \eqref{4.5}) of the well-known  theorem (by Privalov \cite{Pri}) on limiting values of an integral of the form
\begin{equation}f(z)=\int_{a}^{b}\frac{p(t)}{t-z}dt .\label {4.29}\end{equation}
\end{Rk}

Liv\v{s}ic proved \cite{Liv2}   that the multiplicative representation  of the  characteristic matrix function $W(z)$ of the  \eqref{4.5} type gives a method
to construct the corresponding triangular model.  Let us explain this construction.
Putting
\begin{equation}\Phi(\zeta,x)=S_{\zeta}^{-1}F_{1}(x), \label{4.23}\end{equation}
where $S_{\zeta}^{-1}$ is applied to $F_1$ columnwise, we have
\begin{equation}\sigma(\zeta)=iJ\frac{1}{2\pi}\int_{a}^{\zeta}F_{1}(x)\Phi^{*}(\zeta,x)dx=iJ\sigma_{1}(\zeta).
\label{4.24}\end{equation}
\begin{Rk}\label{RkN} The $S$-node considered in this section is symmetric, that is,
$$S_b=S_b^*, \quad B=A^*, \quad \Pi_2^*=iJ\Pi_1^*, \quad \Gamma_2^*=iJ\Gamma_1^*.$$
It follows that
\begin{equation}\sigma_1(\zeta)^*=\sigma_1(\zeta),\label {n6}\end{equation}
and
\begin{equation}W(\overline{z})^*JW(z)=J.\label {n7}\end{equation}
\end{Rk}
We need the following lemma.
\begin{La}\label{Lemma 4.6}  \cite{Sak6} The function $\sigma_{1}(\zeta)$ is
a monotonically increasing matrix function of bounded variation.\end{La}
\emph{We assume that $\sigma_{1}(\zeta)$ is absolutely continuous.}
Thus, $\sigma_{1}^{\prime}(x)$ is well defined and positive definite: $\sigma_{1}^{\prime}(x){\geq}0.$
Introduce a matrix function $\beta(x)$ such that
\begin{equation}\beta^{*}(x)\beta(x)=\sigma_{1}^{\prime}(x).\label{4.25}\end{equation}
\begin{Tm}\label{Theorem 4.7} Let  the conditions \eqref{4.2}--\eqref{4.4} be fulfilled and let $\sigma_{1}(\zeta)$ be absolutely continuous.
Then, we have:\\
1. The matrix function $W(z)$  is the characteristic matrix function of the triangular Liv\v{s}ic model
\begin{equation}Af=xf(x)+i\int_{a}^{x}f(t)\beta(t)Jdt\beta^{*}(x),\quad a{<}x{<}b.\label{4.26}\end{equation}
2. The matrix function $W(z)$ is the  monodromy matrix of the differential system
\begin{equation}\frac {d}{dx}W(x,z)=iJ\frac{\sigma_{1}^{\prime}(x)}{x-z}W(x,z),\quad W(a,z)=I,\label{4.27}\end{equation}
that is,
\begin{equation}W(z)=W(b,z).\label{4.28}\end{equation}\end{Tm}
\begin{Cy}\label{Crollary 4.8} Formulas \eqref{4.6} and \eqref{4.24} give  a procedure to recover  system \eqref{4.27} from the corresponding monodromy matrix
$W(z)=W(b,z)$.\end{Cy}

\begin{Rk}\label{Remark 4.10}The multiplicative representation \eqref{4.1} where $\alpha(x)=0$ was considered in detail in the book \cite{Sak8}.\end{Rk}
\section{ Wiener--Masani prediction theory}
Let us consider the case, where $J=I$ and the relation
\begin{equation} \lim W(z)=I \quad (z{\to}\infty) \label{5.1}\end{equation}
is fulfilled. Using polar decomposition
\eqref{1.4}, we obtain
\begin{equation}W_{+}^{*}(x)W_{+}(x)=R^{2}(x), \label{5.2}\end{equation}
where
\begin{equation}R(x){\geq}I\quad {\mathrm{for}} \quad  x{\in}[a,b], \quad R(x)=I \quad {\mathrm{for}} \quad x{\notin}[a,b].\label{5.3}\end{equation}
In the prediction theory, the following Wiener--Masani problem plays a fundamental role (see \cite{WiMa, Ros}):

\emph{Let $R(x)$ satisfy the inequalities $I{\leq}R(x){\leq}MI$. Find the maximal $W(z)$ such that \eqref{5.2} holds.}\\
Here,  by maximal $W(z)$ we understand $W(z)$ satisfying
\begin{equation}W^{*}(z)W(z)\geq V(z)^*V(z) \,\, {\mathrm{for\,\, all}} \,\, V(z), \,\, {\mathrm{such \,\, that}} \,\, V_{+}^{*}(x)V_{+}(x)=R^{2}(x),
 \label{n5}\end{equation}
and all $z\in \BC_+$. (It is assumed also that the matrix functions $V(z)$ are analytic in $\BC_+$  and $\lim_{z\to \infty}V(z)=I$.)

  N. Wiener and P. Masani proved the existence of the solution $W(z)$ of the formulated problem \cite{WiMa}. However,
their method of constructing  $W(z)$
is not effective. We believe that our methods of constructing $W(z)$ (see section 4) are more simple and  effective. In this section,
it remains to prove that the constructed solution
$W(z)$ is maximal.
\begin{Tm}\label{Theorem 5.1} The matrix function $W(z)$ given by \eqref{4.5} (and \eqref{4.6} or, equivalently, \eqref{4.24}) is maximal.\end{Tm}
\emph{Proof.} According to \eqref{4.2} and \eqref{n7}, $W^{-1}(z)=JW(\overline{z})^*J$ exists and is bounded in the half-plane $\BC_+$.
Taking some other solution $V(z)$ of the Wiener-Masani problem
consider the matrix function $U(z)=V(z)W^{-1}(z)$. Then, we have
\begin{equation}U_{+}^{*}(x)U_{+}(x)=I.\label{5.4}\end{equation}
It follows from \eqref{5.4} and maximum theorem for analytical functions that
\begin{equation}U^{*}(z)U(z){\leq}I,\label{5.5}\end{equation}
which yields the inequality
\begin{equation}W^{*}(z)W(z){\geq}V^{*}(z)V(z).\label{5.6}\end{equation}
The theorem is proved.
\begin{Cy}\label{Corollary 5.2} Let conditions of Theorem \ref{Theorem 5.1}  be fulfilled. Then, the Wiener-Masani problem has only one
(maximal) solution.\end{Cy}
\emph{Proof.}  We showed that $W(z)$ is a solution of the Wiener-Masani problem.
Assume, that $V(z)$ is  another solution of the Wiener-Masani problem. Hence, we have
\begin{equation}W^{*}(z)W(z)=V^{*}(z)V(z).\label{5.7}\end{equation} Relation \eqref{5.7} implies that
\begin{equation}U^{*}(z)U(z)=I.\label{5.8}\end{equation} Taking into account \eqref{5.1} and \eqref{5.8} and using maximum principle, we obtain
$U(z)\equiv I$ (i.e., $V(z)\equiv W(z)$). The corollary is proved.
\begin{Rk}\label{Remark 5.3}The Wiener-Masani problem was solved in  \cite{WiMa} without the assumption that $W(z)$ is  regular at $z=\infty.$ In our approach, one can  omit this assumption as well (see \cite{Sak7}).\end{Rk}
Using  Corollary 5.2 we derive the next assertion.
\begin{Tm}\label{Theorem 5.4} Suppose that  a matrix function $W(z)$ satisfies
\eqref{4.2} and admits two representations
\begin{equation}W(z)=\int \limits_a^{\overset{b}{\curvearrowleft}} \, \exp\left[ -i\frac{d\sigma_{k}(t)}{t-z}\right]\quad (k=1,2),\label{5.9}\end{equation}
where $\sigma_{k}(t)$ are continuous monotonically increasing matrix functions of bounded variation with $\sigma_{1}(a)=\sigma_{2}(a).$ Then,
\begin{equation} \sigma_{1}(t)=\sigma_{2}(t),\quad a{\leq}t{\leq}b.\label{5.10}\end{equation}\end{Tm}
Theorem  \ref{Theorem 5.4} shows the uniqueness of the system \eqref{4.27} with  a given monodromy matrix $W(z)$. Yu.P. Ginzburg earlier
proved  \cite{Gin} Theorem \ref{Theorem 5.4}  using other methods.
\section{A special class of $J$-modules \\ and random matrices.}
V.P. Potapov \cite{Pot} introduced and investigated the $J$-modules  $R$ of $m{\times}m$ matrices.
Recall that the $J$-module $R$ of an $m{\times}m$ matrix  has the following properties:

 1. The eigenvalues   of $R$ are positive.

 2. The matrix $R$ satisfies the relations
 \begin{equation} JR=R^{*}J,\quad RJR^{*}{\geq}J.\label{6.1}\end{equation}
 where $J=J^{*}$ and $J^{2}=I$.

 In this section, we shall consider a special Z-class (see Definition \ref{Definition 1.2}) of $J$-modules $R$.
Namely,  we assume that $R$ satisfy the additional condition
 \begin{equation}(R-I)^{2}=0,\label{6.2}\end{equation}
 and that
 \begin{equation}
 m=2, \quad
 J=\left(
                     \begin{array}{cc}
                       -1 & 0 \\
                       0& 1 \\
                     \end{array}
                   \right)\label{6.3}.\end{equation}
 It follows from \eqref{6.1}--\eqref{6.3}, that $R$  has the form
 \begin{equation}R=\left(
                     \begin{array}{cc}
                       1-|\psi| & \psi \\
                       -\overline{\psi}& 1+|\psi| \\
                     \end{array}
                   \right)\label{6.4}\end{equation}
 We note that if $R$ belongs to the Z-class, then $R^{2}$ belongs to the Z-class as well. Indeed, \eqref{6.4} implies that
\begin{equation}R^{2}=\left(
                     \begin{array}{cc}
                       1-2|\psi| & 2\psi \\
                       -2\overline{\psi}& 1+2|\psi| \\
                     \end{array}
                   \right).\label{6.5}\end{equation}
Further we need the following relations:
\begin{equation}R^{-1}=\left(
                     \begin{array}{cc}
                       1+|\psi| & -\psi \\
                       \overline{\psi}& 1-|\psi| \\
                     \end{array}
                   \right)\label{6.6}\end{equation}
and
 \begin{equation}J[R-R^{-1}]=2\left(
                     \begin{array}{cc}
                       |\psi| & -\psi \\
                       -\overline{\psi}& |\psi| \\
                     \end{array}
                   \right).\label{6.7}\end{equation}
Let us introduce the vector function
\begin{equation}F_{1}(x)=\col[\phi(x),\, -\overline{\phi(x)}],\quad {\mathrm{where}} \quad \phi^{2}(x)=2\psi(x). \label{6.8}\end{equation}
It is easy to see that the vector  $F_{1}(x)$ satisfies the following relations:
\begin{equation}F_{1}(x)F_{1}^{*}(x)=J[R(x)-R^{-1}(x)],\quad  F_{1}^{*}(x)JF_{1}(x)=0.\label{6.9}
\end{equation}
Now, we  can solve the following\\
\emph{Problem I.}\\
\emph{Let $J$-module $R(x)$ be  bounded and  satisfy  \eqref{6.2} and \eqref{6.3} (i.e., let the bounded $R$ belong to Z-class). Construct the corresponding matrix
function $W(z)$.}

Using \eqref{4.12}, \eqref{4.20} and \eqref{6.9}, we obtain
\begin{equation}S_{\zeta}f=f(x)+\frac{i}{2\pi}\dashint_a^\zeta\frac{\phi(x)\overline{\phi(t)}-
\overline{\phi(x)}\phi(t)}{x-t}f(t)dt.
\label{6.10} \end{equation}
Our considerations above yield the following assertions.
\begin{Pn}\label{Proposition 6.1}If  $\psi(x)$ is bounded,   then the matrix functions
$R(x)$, $R^{-1}(x)$ and $F_{1}(x)$ are bounded as well.\end{Pn}
\begin{Pn}\label{Proposition 6.2}
If  $\psi(x)$ is bounded,   then the operator $S_{\zeta}$ is bounded in the space $L^{2}(a,\zeta)$ $\,(\zeta{\in}(a,b])$.\end{Pn}
\begin{Tm}\label{Theorem 6.3}If  $\psi(x)$ is bounded  and the operator
$S_{b}$   is positive and invertible, then the solution $W(z)$ of Problem I has the form
\begin{equation}
W(z)=\int \limits_a^{\overset{b}{\curvearrowleft}} \, \exp\left[-iJ \frac{d\sigma_{1}(t)}{t-z}\right],\label{6.11}
\end{equation}
where
\begin{equation}\sigma_{1}(\zeta)=\frac{1}{2\pi}\int_{a}^{\zeta}F_{1}(x)\Phi^{*}(\zeta,x)dx
\label{6.12}\end{equation}
and
\begin{equation}\Phi(\zeta,x)=\col[\Phi_{1}(\zeta,x),\, -\overline{\Phi_{1}(\zeta,x)}]=S_{\zeta}^{-1}F_{1}(x). \label{6.13}\end{equation}\end{Tm}
Formula \eqref{6.12} may be rewritten in the form
\begin{equation}\sigma_{1}(\zeta)=\frac{1}{2\pi}\int_{a}^{\zeta}
\begin{pmatrix}
    \phi(x)\overline{\Phi_{1}(\zeta,x)} &\,\,-\phi(x)\Phi_{1}(\zeta,x)  \\ \\
  -  \overline{\phi(x)\Phi_{1}(\zeta,x)} & \,\,  \overline{\phi(x)}\Phi_{1}(\zeta,x)
   \end{pmatrix}
 dx.
\label{6.14}\end{equation}
\begin{Rk}\label{Remark 6.4} Problem I is equivalent to Riemann-Hilbert problem \eqref{4.22}.
\end{Rk}
Now, suppose that  $S_{b}$ admits  triangular factorization
\cite{GoKr4}, Ch.IY and \cite{Sak12}, that is, $S_b$ may be represented in the form
\begin{equation} S_b=S_{-}^{*}S_{-}, \label{6.15}\end{equation}
where $S_{-}^{\pm1}$ are bounded operators and
\begin{equation}S_{-}^{\pm1}\clp_{\mu}=\clp_{\mu}S_{-}^{\pm1}\clp_{\mu}.\label{6.16}\end{equation}
Here, $\big(\clp_{\mu}f\big)(x)=f(x)\,$ for $\,\mu{\leq}x{\leq}b\,$ and  $\,\big(\clp_{\mu}f\big)(x)=0\,$ for $\,a{\leq}x{\leq}\mu$,\\
$f(x){\in}L^{2}(a,b).$
From  \eqref{6.14}---\eqref{6.16} we deduce the next assertion.
\begin{Tm}\label{Theorem  6.5} If the conditions of Theorem \ref{Theorem  6.3} are fulfilled and the corresponding operator $S_{b}$ admits the triangular factorization, then
\begin{equation}\sigma_{1}(\zeta)=\frac{1}{2\pi}\int_{a}^{\zeta}
\left(
  \begin{array}{cc}
    |q(x)|^2 &  \,\,-q^{2}(x) \\ \\
 -  \overline{q^{2}(x)}  & \,\, |q(x)|^2
  \end{array}
\right)dx,
\label{6.17}\end{equation}
where
\begin{equation}q(x)=S_{-}^{-1}\phi. \label{6.18}\end{equation}\end{Tm}
\emph{Proof.} The kernel of the operator $S_b$ is real-valued (see \eqref{6.10}). Therefore, we can choose
the operator $S_{-}^{-1}$ such that its kernel   is real-valued (real) as well. Hence, we have
\begin{equation}\overline{S_{-}^{-1}\overline{f}}=S_{-}^{-1}f.\label{6.19}\end{equation}
Relations \eqref{6.14}--\eqref{6.16} and \eqref{6.19} yield
\eqref{6.17} and \eqref{6.18}. The theorem is proved.

\begin{Cy}\label{Corollary 6.6} Let conditions of Theorem \ref{Theorem  6.5} be fulfilled. Then, the
matrix function $\sigma_{1}(x)$ is absolutely continuous and
\begin{equation}\sigma_{1}^{\prime}(x)=\frac{1}{2\pi}
\left(
  \begin{array}{cc}
    |q(x)|^2 & \,\,-q^{2}(x) \\ \\
   -\overline{q^{2}(x)}  &\,\, |q(x)|^2\\
  \end{array}
\right).
\label{6.20}\end{equation}\end{Cy}
Taking into account \eqref{6.20}, we obtain  the following assertion.
\begin{Cy}\label{Corollary 6.7} Let the conditions of Theorem \ref{Theorem  6.5} be fulfilled. Then,
the equality
\begin{equation}[J\sigma_{1}^{\prime}(x)]^{2}=0 \label{6.21}\end{equation}
is valid.\end{Cy}
In our paper \cite{Sak17}, we found \cite[Proposition 4.4]{Sak17} the conditions under which
\eqref{6.2} follows from \eqref{6.21}. In the same paper,  we formulated  an {\it open problem}:\\

Find the conditions under which formula \eqref{6.21} follows from \eqref{6.2}.\\
\emph{Corollary 6.7 gives the solution of this problem.}

\begin{Tm}\label{Theorem 6.8} Let the inequality
\begin{equation}\int_{a}^{b}\int_{a}^{b}\left|\frac{\phi(x)\overline{\phi(t)}-
\overline{\phi(x)}\phi(t)}{x-t}\right|^{2}dtdx<\infty \label{6.22}\end{equation}
hold.
Then, the operators $S_{\zeta}$  $\,(a<\zeta{\leq}b)$ of the form \eqref{6.10} admit the triangular factorization \eqref{6.15}.\end{Tm}
This theorem follows directly from the factorization theorem by M.G. Krein    \cite[Chapter 4]{GoKr4}.
\begin{Cy}\label{Corollary 6.9} If  \eqref{6.22} is valid, then there
exists $r$ $(a<r\leq b)$ such that the operators $S_{\zeta}$ of the form \eqref{6.10} are positive, invertible in the space $L^{2}(a,\zeta)$ and  factorable for all $a<\zeta{\leq}r$.\end{Cy}

The multiplicative representation
 \begin{equation}
W(x,z)=\int \limits_a^{\overset{x}{\curvearrowleft}} \, \exp\left[ -\frac{iJd\sigma_{1}(t)}{t-z}\right],\quad z{\notin}[a,b] \label{6.23}
\end{equation}
yields an asymptotic formula
 \begin{equation}W(x,z)=I+M(x)/z+O(1/{z^2}),\quad z{\to}\infty,
\label{6.24}\end{equation}
where
\begin{equation}M(x)=i\int_{a}^{x}Jd\sigma_{1}(t).\label{6.25}\end{equation}
\begin{Rk}\label{Remark 6.10} The matrix function $M(x) $plays an essential role in the random matrix theory  \cite{DIZ}.\end{Rk}
If we represent $\phi(x)$ as the sum of its real and imaginary parts
\begin{equation} \phi(x)=A(x)+iB(x),\quad A(x)=\overline{A(x)},\quad B(x)=\overline{B(x)},
\label{6.26}\end{equation}
formula \eqref{6.10} for the operator $S_{\zeta}$ may be rewritten in a convenient form
\begin{equation}S_{\zeta}f=f(x)-\frac{1}{\pi}\dashint_a^\zeta\frac{A(x)B(t)-
B(x)A(t)}{x-t}f(t)dt.
\label{6.27} \end{equation}

Now, we shall consider some important examples, where condition \eqref{6.2} is fulfilled. Then,  $R(x)$ and $R^{2}(x)$ have  the forms \eqref{6.4} and \eqref{6.5}, respectively.
\begin{Ee}\label{Example 6.11}{Sine-kernel.}\end{Ee}
Consider the case where
\begin{equation} \psi(x )=\frac{1}{2}\gamma{e^{2ix }}, \quad 1\geq{\gamma}>0.\label{6.28}\end{equation}
It follows from \eqref{6.5} and \eqref{6.28} that
\begin{equation}R^{2}(x )=\left(
                           \begin{array}{cc}
                             1-\gamma & \gamma{e^{2ix }} \\
                             -\gamma{e^{-2ix }} & 1+\gamma \\
                           \end{array}
                         \right)\label{6.29}
                         \end{equation}
The Riemann--Hilbert problem \eqref{4.22}, \eqref{6.29} plays an essential role in the random matrix theory \cite{DIZ}.
In view of \eqref{6.8} and \eqref{6.28}, we have
\begin{equation} \phi(x )=i\sqrt{\gamma}\, e^{ix }.\label{6.30}\end{equation}
According to \eqref{6.10} and \eqref{6.30},
the operator $S_{\zeta}$ $($acting in $L^{2}(a,\zeta))$ has the form
 \begin{equation}S_{\zeta}f=f(x)-\gamma\int_{a}^{\zeta}\frac{\sin(x-t)}{\pi(x-t)}f(t)dt. \label{6.31}\end{equation}
\begin{Pn}\label{Proposition 6.12}The operator $S_{\zeta}$, defined by \eqref{6.31}, is positive, invertible, and factorable, and the
 operator $K_{\zeta}=I-S_{\zeta}$ belongs to the trace class.\end{Pn}
\emph{Proof.} It is well known that the operator $S_{\zeta}$ is positive, invertible and that $K_{\zeta}$ belongs to the trace class \cite[p. 167]{DIZ}.
In the case \eqref{6.29}, the condition \eqref{6.22} is fulfilled. Hence, the operator $S_{\zeta}$ is factorable. The proposition is proved.
\begin{Cy}\label{Corollary 6.13} Let \eqref{6.28} be valid. Then, $W(z)$ may be recovered using \eqref{6.13}, \eqref{6.14} $($or
\eqref{6.17}, \eqref{6.18}$)$ and \eqref{6.23}. Its asymptotics is given by
\eqref{6.24} and \eqref{6.25}.
\end{Cy}
\begin{Rk}\label{Remark 14}
Let $U=\bigcup_{k=1}^{n}(a_k,b_k)$
 be a  union of n disjoint intervals in $\mathbb{R}$. The results and the methods of section 6 may be also used when we consider the domain U.\end{Rk}
\begin{Ee}\label{Example 6.15} {Bessel kernel}. \end{Ee}
Consider the case
\begin{equation}\phi(x )=A(x )+iB(x ), \quad A(x )=\sqrt{\gamma\pi}J_{\alpha}(x ),\quad B(x)=xA^{\prime}(x ),\label{6.32}\end{equation}
where $x{\geq}0$ and $J_{\alpha}(x)$ is the Bessel function of the order $\alpha,\, \alpha>-1/2$.
Thus, we have the Riemann-Hilbert problem \eqref{4.22}, where the matrix $R^{2}(x )$ is determined by \eqref{6.5}, \eqref{6.32}  and the relation
\begin{equation}\psi(x )=\frac{1}{2}\phi^{2}(x ).\label{6.34}\end{equation}
 Let us write the corresponding operator $S_{\zeta}$:
\begin{equation}S_{\zeta}f=f(x)+\int_{0}^{\zeta}k(x,t)f(t)dt,\,f(t){\in}L^{2}(0,\zeta),\quad
0<\zeta<\infty,\label{6.35}\end{equation}
where
\begin{equation}k(x,t)=\frac{i}{2\pi}\frac{\phi(x)\overline{\phi(t)}-
\overline{\phi(x)}\phi(t)}{x-t}=-\frac{A(x)B(t)-A(t)B(x)}{x-t}.
\label{6.36} \end{equation}
\begin{Pn}\label{Proposition 6.16}The operator $S_{\zeta}$, defined by \eqref{6.32}, \eqref{6.35} and \eqref{6.36} is positive, invertible and
factorable,
 and the operator $K_{\zeta}=I-S_{\zeta}$ belongs to the trace class.\end{Pn}
Proof. Taking into account Christoffel-Darboux formula, we have \cite{TrWi}
\begin{equation}k(x,t)=-\frac{\gamma}{4}\int_{0}^{1}J_{\alpha}(\sqrt{xs})J_{\alpha}(\sqrt{ts})ds.
\label{6.37}\end{equation}
The
 operator $K_{\zeta}=I-S_{\zeta}$ is self-adjoint. It follows also from \eqref{6.35} and \eqref{6.37} that
 \begin{equation}(K_{\zeta}f,f)=\frac{\gamma}{4}\int_{0}^{1}\left|\int_{0}^{\zeta}J_{\alpha}(\sqrt{xs})f(x)dx\right|^{2}ds.
\label{6.38}\end{equation} The last relation may be rewritten in the form
\begin{equation}(K_{\zeta}f,f)=2\gamma\int_{0}^{1}\left|\int_{0}^{\sqrt{\zeta}}J_{\alpha}(t\eta)f(t^2)tdt\right|^{2}\eta d\eta ,
\label{6.39}\end{equation}
where $\sqrt{x}=t,\,\sqrt{s}=\eta .$ We introduce the  Hankel transformation
\begin{equation}F_{\alpha}(\eta )=\int_{0}^{\sqrt{\zeta}}J_{\alpha}(t\eta )f(t^2)tdt .\label{6.40}\end{equation}
According to \eqref{6.39} and \eqref{6.40}, we have
\begin{equation}(K_{\zeta}f,f)=2\gamma\int_{0}^{1}\left|F_{\alpha}(\eta )\right|^{2}\eta d\eta ,
\label{6.41}\end{equation}
Hence, the operator $K_{\zeta}$ is non-negative. It follows from \eqref{6.37}
 that
\begin{equation}\int_{0}^{\zeta}k(x,x)dx<\infty.\label{6.42}\end{equation}
Thus, the operator $K_{\zeta}$  belongs to the trace class.
The Hankel transformation is unitary.  Thus, we derive
\begin{equation}(K_{\zeta}f,f){\leq}2\gamma\int_{0}^{\infty}|F_{\alpha}(\eta )|^{2}\eta d\eta =
2\gamma\int_{0}^{\sqrt{\zeta}}|f(t^2)|^2tdt.
\label{6.43}\end{equation}
Hence, we obtain the inequality
\begin{equation}(K_{\zeta}f,f){\leq}\gamma\int_{0}^{\zeta}|f(x)|^2dx.\label{6.44}\end{equation}
Relation \eqref{6.44} implies that
\begin{equation} \| K_{\zeta}\| {\leq}\gamma.\label{6.45}\end{equation}
  It follows from \eqref{6.37}
   that condition \eqref{6.22} is fulfilled. If $0<\gamma<1$ then the corresponding operator $S_{\zeta}$ is positive,
invertible and factorable.

Let us consider separately the case, where $\gamma=1$. Assume that $\| K_{\zeta}\| =1$. In this case, for some $f{\ne}0$  we have the equality
$K_{\zeta}f=f$. From  \eqref{6.41} and \eqref{6.43} we deduce that $F_{\alpha}(\eta )=0$ if $\eta >1$.  Taking into account the analyticity of
 $F_{\alpha}(\eta )$  in the domain
 $\Re(\eta )>0$, we see that
$F_{\alpha}(\eta )=0$ for $\eta >0$. It means that $f=0$, and we arrive at a contradiction. Thus, $\| K_{\zeta}\| <1$.
The proposition is proved.
\begin{Cy}\label{Corollary 6.17} Let \eqref{6.32} be valid. Then, $W(z)$ may be recovered using \eqref{6.13}, \eqref{6.14} $($or
\eqref{6.17}, \eqref{6.18}$)$ and \eqref{6.23}. Its asymptotics is given by
\eqref{6.24} and \eqref{6.25}.
\end{Cy}
\begin{Ee}\label{Example 6.18}Airy kernel.\end{Ee}
Let us consider the operator
\begin{equation}S_{\zeta}f=f(x)-\int_{\zeta}^{\infty}k(x,t)f(t)dt,\quad
f(x){\in}L^{2}(\zeta,\infty),\quad \zeta >0, \label{6.46}\end{equation}
where the  kernel $k(x,t)$ has the form
\begin{equation} k(x,t)=\frac{A_{i}(x)A^{\prime}_{i}(t)-A^{\prime}_{i}(x)A_{i}(t)}{x-t}.\label{6.47}\end{equation}
Here, $A_{i}(x)$ is Airy function. Since  $S_{\zeta}$ is acting in $L^{2}(\zeta,\infty)$ instead of $L^{2}(a,\zeta)$
(as earlier), the projectors in the analog of \eqref{4.17+} for this case differ from $P_{\zeta}$.

The following estimates  of Airy functions are valid (see \cite[Ch. 11]{Olv}):
\begin{equation}A_{i}(x)=O\big(x^{-1/4}\exp[(-2/3)x^{3/2}]\big),\quad A^{\prime}_{i}(x)=
O\big(x^{1/4}\exp[(-2/3)x^{3/2}]\big),\label{6.48}\end{equation}
where $x{\to}\infty,\,\, |\arg({x})|<\pi.$ We take also into account  that
\begin{equation}|A_{i}(0)|+|A^{\prime}_{i}(0)|<\infty.\label{6.49}\end{equation}
We shall need  the Christoffel-Darboux type formula (see
\cite{TrWi}):
\begin{equation}k(x,t)=\int_{0}^{\infty}A_{i}(x+s)A_{i}(t+s)ds.\label{6.50}\end{equation}
\begin{Pn}\label{Proposition 6.19} The operator $S_{\zeta}$, defined by \eqref{6.46},  is positive, invertible and factorable, and the operator $K_{\zeta}=I-S_{\zeta}$
 belongs to the trace class.\end{Pn}
\emph{Proof.}
The operator $K_{\zeta}$ is self-adjoint. Using \eqref{6.50} we have
\begin{equation}(K_{\zeta}f,f)_{\zeta}=\int_{0}^{\infty}|F(s)|^{2}sds,
\label{6.51}\end{equation}where
\begin{equation}F(s)=\int_{\zeta}^{\infty}A_{i}(x+s)f(x)dx, \label{6.52}\end{equation}
and the symbol $(f,g)_{\zeta}$ denotes the inner product in the space $L^{2}(\zeta,\infty)$.
Hence, the operator $K_{\zeta}$ is non-negative. It follows from \eqref{6.48}
and \eqref{6.50} that
\begin{equation}\int_{\zeta}^{\infty}k(x,x)dx<\infty.\label{6.53}\end{equation}
Thus, the operator $K_{\zeta}$  belongs to the trace class.
Since the  transformation \eqref{6.52} is unitary, we have
\begin{equation}(K_{\zeta}f,f){\leq}\int_{-\infty}^{\infty}|F(s)|^{2}ds=
\int_{\zeta}^{\infty}|f(x)|^2dx.
\label{6.54}\end{equation}
Hence, $\| K_{\zeta}\| {\leq}1.$ It is easy to see that the analog of \eqref{6.22} is fulfilled:
\begin{equation}\int_{0}^{\infty}\int_{0}^{\infty}\left|\frac{\phi(x)\overline{\phi(t)}-
\overline{\phi(x)}\phi(t)}{x-t}\right|^{2}dtdx<\infty .\label{6.22'}\end{equation}
If $\| K_{\zeta}\| <1$,  then the corresponding operator $S_{\zeta}$ is positive,
invertible and factorable.

Now, let us assume that $\| K_{\zeta}\| =1.$  In this case,  for some $f{\ne}0$ we have
$K_{\zeta}f=f$. From  \eqref{6.51} and \eqref{6.54} we deduce that $F(s)=0$ if $s<0$. Since the function $F(s)$ is analytic in the domain $-\zeta<\Re({s})<\infty$ (see \eqref{6.52}),
the equality $F(s)=0$ for $s \geq 0$ holds. It means that $f=0$. We arrived at the contradiction. Thus, $\| K_{\zeta}\| <1$. The proposition is proved.

\section{On the connection between Fredholm \\ determinant and Riemann--Hilbert problem}
In this section, we assume that the conditions of Theorem \ref{Theorem 6.3} are fulfilled and that the operator $(I-S_{\zeta})$ $(a<\zeta<b)$,
where  $S_{\zeta}$ is given by \eqref{6.10}, belongs to the trace class.
Then, the following assertion is valid.
\begin{Pn}\label{Proposition 7.1}
The operator $S_b$ is factorable and
 Fredholm determinant ${\bf P}_{\zeta}=\det\big(S_{\zeta}\big)$ exists.\end{Pn}
We introduce the operator
\begin{equation}Q_{\zeta}f=\int_{0}^{1}k(x\zeta,t\zeta){f(t)}d(t\zeta),
\label{7.1}\end{equation}
where the kernel $k(x,t)$ is defined by \eqref{6.35}. It is easy to see that the operators
$S_{\zeta}$ and $I+Q_{\zeta}$ have the same eigenvalues. Hence,
\begin{equation}{\bf P}_{\zeta}=\det(I+Q_{\zeta}).\label{7.2}\end{equation}
Using the equality $\log{{\bf P}_{\zeta}}=\tr[\log(I+Q_{\zeta})]$, we derive
\begin{equation}\frac{d}{d\zeta}\log{{\bf P}_{\zeta}}=
\tr[(I+Q_{\zeta})^{-1}\frac{d}{d\zeta}Q_{\zeta}].\label{7.3}
\end{equation}
\begin{Ee}\label{Example 7.2}Sine kernel.\end{Ee}
Let us consider again the case where
\begin{equation}\phi(x\zeta)=i\sqrt{\gamma}e^{ix\zeta}.\label{7.4}\end{equation}
 Taking into account \eqref{6.35} and \eqref{7.1}, we have
 \begin{equation}Q_{\zeta}f=-\gamma\int_{0}^{1}\frac{\sin\zeta(x-t)}{\pi(x-t)}f(t)dt, \label{7.5}\end{equation}
 where the operators on both sides of \eqref{7.5} act in $L^{2}(0,1)$. Relations \eqref{6.8} and
\eqref{7.4}, \eqref{7.5} imply that
\begin{equation}\frac{d}{d\zeta}Q_{\zeta}f=-\frac{\gamma}{2\pi}\int_{0}^{1}
F_{1}^{*}(t\zeta)F_{1}(x\zeta)f(t)dt.\label{7.6}\end{equation}
It follows from \eqref{7.1} and \eqref{7.6} that
\begin{equation} \tr[(I+Q_{\zeta})^{-1}\frac{d}{d\zeta}Q_{\zeta}]=
-\frac{\gamma}{2\pi}\int_{0}^{1}F_{1}^{*}(t\zeta)\Phi(t\zeta,\zeta)d(t\zeta),\label{7.7}\end{equation}
where
\begin{equation}\Phi(x\zeta,\zeta)=(I+Q_{\zeta})^{-1}F_{1}(x\zeta).\label{7.8}
\end{equation}
From \eqref{6.12}, \eqref{7.7} and \eqref{7.8} we derive
\begin{equation}\tr[(I+Q_{\zeta})^{-1}F_{1}(x\zeta)]=-\gamma\sigma_{1}(\zeta).\label{7.9}
\end{equation}
Finally, using equalities  \eqref{7.3} and (7.9) we obtain the well-known result \cite{DIZ}:
\begin{equation}\frac{d}{d\zeta}\log{{\bf P}_{\zeta}}=i[m_{2,2}(\zeta)-m_{1,1}(\zeta)],
\label{7.10}\end{equation} where $m_{i,j}(\zeta)$ are the elements of the matrix $M(\zeta)$ given by \eqref{6.25}.

\section{Triangular form}
{\bf 1.} During the final year of my postgraduate studies (in 1956), I read the Aronszajn--Smith paper \cite{AS}. The paper contained the following important result:\\
\emph{Each compact operator in a Hilbert space has a nontrivial invariant subspace.}

Together with my wife E. Melnichenko, we translated this paper into Russian and published this translation in the Russian journal ``Mathematics".
In the same 1956 year, M. Liv\v{s}ic gave me his manuscript, where he independently proved the Aronszajn--Smith theorem.  These two proofs were very similar, although  M. Liv\v{s}ic used the notion of the characteristic matrix function and N. Aronszajn and K.T. Smith did not.
Taking into account these works, I decided  to use geometrical methods (in particular, Aronszajn--Smith theorem) for solving the problems connected with triangular forms \cite{Sak15, Sak16}. More precisely, for the results in this paragraph see \cite{Sak15} and for results in the next paragraph {\bf 2} see \cite{Sak16}.

Let us consider linear operators $A$ acting in a separable Hilbert space $H$.
\begin{Dn}\label{Definition 8.1}
The maximal invariant subspace of  the  operators $A$ and $A^{*}$, on which the equality
$AA^{*}=A^{*}A$ is fulfilled,  is called the additional component.\end{Dn}
\begin{Dn}\label{Definition 8.2} We say that the operator $A$ $($acting in a separable Hilbert space $H)$
belongs to the class $T$ if for any two invariant subspaces $H_1$  and  $H_2$  $(H_1{\subseteq}H_2$,
$\dim(H_2{\ominus}H_1)>1)$ of the operator
$A$ there exists a third invariant subspace  $H_3$ of $A$ such that
$H_1{\subseteq}H_3{\subseteq}H_2$
and  $H_1{\ne}H_3$, $H_3{\ne}H_2$.\end{Dn}
It follows  from the Aronszajn--Smith theorem \cite{AS} that all compact operators belong
to the class T.
\begin{Dn}\label{Definition 8.3} We say that the operator $\vec{A}\in [L^2_m(0,1),L^2_m(0,1)]$ is triangular  if it has the following form
\begin{equation}\vec{A}f=\frac{d}{dx}\int_{0}^{x}f(t)K(x,t)dt,\quad 0{\leq}x{\leq}1,
\label{8.1}\end{equation}
where
\begin{equation} f(t)=[f_1(t),f_2(t),...,f_m(t)],\quad m{\leq}\infty.\label{8.2}\end{equation}\end{Dn}

We denote the additional component of $A$ by ${\mathcal{H}}_A$. Clearly, $H{\ominus}{\mathcal{H}}_A$ is an invariant subspace of $A$.
\begin{Tm}\label{Theorem 8.4} For every operator $A$ from the class $T$, there exist a triangular operator $\vec{A}$ of the form \eqref{8.1} and a unitary operator $U$ with the
following properties:\\
1. The operator $U$ is a one to one mapping from the space $H{\ominus}{\mathcal{H}}_A$ onto the space
$L^2_m(0,1){\ominus}{\mathcal{H}}_{\vec{A}}$, where ${\mathcal{H}}_A$ and ${\mathcal{H}}_{\vec{A}}$ are the additional components of
$A$ and $\vec{A}$, respectively.\\
2. The operator ${\mathcal{A}}$, which is a reduction of $A$ onto $H{\ominus}{\mathcal{H}}_A$, is transformed by $U$ into the operator $\vec{\mathcal{A}}=U{\mathcal{A}}
U^{-1}$ which is a reduction of $\vec{A}$ onto $L^2_m(0,1){\ominus}{\mathcal{H}}_{\vec{A}}.$
\end{Tm}
\begin{Cy}\label{Corollart 8.6} Any compact operator is unitary equivalent $($up to additional component$)$ to the operator of the form \eqref{8.1}.\end{Cy}
{\bf 2.} Further in the text, we  consider operators of the form
\begin{equation}A=L+iB, \label{8.3}\end{equation}
where $L$ is a bounded self-adjoint operator and $B$ is a  self-adjoint operator from the Hilbert--Schmidt class.
\begin{Dn}\label{Definition 8.7} The class of operators $A$  of the form \eqref{8.3} is denoted by $T_{1}$.\end{Dn}
\begin{Tm}\label{Theorem 8.8}  Every operator of the class $T_1$  has a non-trivial invariant subspace .\end{Tm}
\begin{Cy}\label{Corollary 8.9}The class $T_1$  is contained in the class $T$ $(T_1\subset T)$.\end{Cy}
Note that the class $i\Omega$ (of the operators $iA$, where $A$ has the form \eqref{8.3} with $B$ belonging to the trace class) is contained in the class $T_1$ (see \cite{Liv2}).

If $A\in T_1$, its triangular model, that is, the operator $\vec{A}$ which appears in Theorem \ref{Theorem 8.4} may be described in greater detail.
\begin{Tm}\label{Theorem 8.10} The triangular model of an operator  $A\in T_1$,
whose spectrum is purely real, can be represented in the form:
\begin{equation}\vec{A}\phi=\phi(x)L(x)+\int_{0}^{x}\phi(t)N(x,t)dt \quad (\phi\in L^2_m(0,1)), \label{8.4}\end{equation}
where $L(x)$ and $N(x,t)$ are $m{\times}m$ matrices $\big(N=\{n_{i,j}\}_{i,j=1}^m\big)$ and
\begin{equation}L(x)=L^{*}(x),\quad \sum_{i,j=1}^{m}\int_{0}^{1}\int_{0}^{1}|n_{i,j}(x,t)|^{2}dt<\infty.\label{8.5}
\end{equation}\end{Tm}
\begin{Tm}\label{Theorem 8.11}If the spectrum of the operator $A\in T_1$ consists only of the zero, then the triangular model of  $A$ has the form
\begin{equation}\vec{A}\phi=\int_{0}^{x}\phi(t)N(x,t)dt, \label{8.4+}\end{equation}\end{Tm}
\begin{Cy}\label{Corollary 8.12} If the spectrum of the operator $A\in T_1$ consists only of  zero, then the operator $A$ belongs to the Hilbert--Schmidt class.\end{Cy}
Corollary \ref{Corollary 8.12} can be rewritten in the following form:
\begin{Cy}\label{Corollary 8.13} If the spectrum of the operator $A$ consists only of zero  and $\Im{(A)}$ belongs to the Hilbert--Schmidt class, then $\Re{(A)}$ belongs to the Hilbert--Schmidt class as well.\end{Cy}
\begin{Rk}\label{Remark 14+} Unfortunately, Liv\v{s}ic's opinion on the above-mentioned results from
section 8 was sharply negative.\end{Rk}
{\bf 3.} Next, let us describe some developments of the results from the previous paragraphs.
\begin{Dn}\label{Definition 8.15} The class of compact operators $A$ such that
\begin{equation}\tr[(A^{*}A)^{p/2}]<\infty \label{8.7}\end{equation}
is denoted by ${\mathfrak{S}}_p$ $(1{\leq}p<\infty)$.
\end{Dn}
I.C. Gohberg and M.G. Krein in \cite{GoKr1} and V.I. Macaev in \cite{Mac} generalized our Corollary \ref{Corollary 8.13}, where the case
$p=2$ is dealt with. The following generalization
is valid.
\begin{Cy}\label{Corollary 8.16} If the spectrum of the operator $A$ consists only of  zero  and $\Im{(A)}$ belongs to  ${\mathfrak{S}}_p\,\,(1{\leq}p<\infty)$, then $\Re{(A)}$ belongs to the ${\mathfrak{S}}_p$ class as well.\end{Cy}
I.C. Gohberg and M.G.Krein \cite{GoKr2}, and V.I. Macaev \cite{Mac} generalized also our Theorem \ref{Theorem 8.8} and Corollary \ref{Corollary 8.9}.
\begin{Tm}\label{Theorem 8.17} If  $A$ is a bounded operator and $\Im{(A)}$ belongs to  ${\mathfrak{S}}_p$\\
 $(1{\leq}p<\infty)$, then the operator $A$ has a non-trivial invariant subspace and belongs to the class $T$.\end{Tm}
We  proved this theorem  for the case $1{\leq}p{\leq}2$.
Later, V.I.  Macaev obtained (in \cite{Mac}) a more general result  than Theorem \ref{Theorem 8.17}.
\begin{Rk}\label{Remark 8.18}  A new interpretation of the triangular form  was found by  M.~Brodskii  in \cite{Bro}.\end{Rk}
\begin{Rk}\label{Remark 8.19} Some ideas from the operator  spectral
theory \cite{Dun} and our triangular representation theory \cite{Sak15, Sak16} were used recently in our paper \cite{Sak18}.\end{Rk}

\section{Reduction to the simplest form}
When I was still a graduate student, M. Liv\v{s}ic  gave me a hand typed copy of his manuscript \cite{Liv1} to write
in the appropriate  formulas by hand (as was done at that time).
I was greatly impressed by that work. I had already a good knowledge of linear algebra and
understood that there was a two-step approach to the corresponding problems:\\
1. Reduction of the matrix to the triangular form.\\
2. Reduction of the matrix to the simplest form.\\
In my opinion,  Liv\v{s}ic made (in his article) the first step, and I decided to try to do the second step.
I started with a simple but non-trivial  example of a triangular Liv\v{s}ic model:
\begin{equation}Af=xf(x)+i\alpha\int_{0}^{x}f(t)dt,\quad 0{<}x{<}\ell, \label{9.1}\end{equation}
where $\alpha$ is real and $f(x){\in}L^{2}(0,\ell)$.
During my final year at the  institute, I proved the following result \cite{Sak10}:
\begin{Tm}\label{Theorem 9.1} The operator $A$ admits  representation of the form
\begin{equation}A=B^{-1}QB, \quad Qf=xf(x),\label{9.2}\end{equation}
where  $B$ and $B^{-1}$  are bounded in $L^{2}(0,\ell)$  operators given by the formulas
\begin{equation} Bf=\frac{1}{\Gamma(i\alpha+1)}\frac{d}{dx}\int_{0}^{x}f(t)(x-t)^{i\alpha}dt,\quad
0{<}x{<}\ell,\label{9.3}\end{equation}
\begin{equation} B^{-1}f=\frac{1}{\Gamma(-i\alpha+1)}\frac{d}{dx}\int_{0}^{x}f(t)(x-t)^{-i\alpha}dt,\quad
0{<}x{<}\ell.\label{9.4}\end{equation}
Here, $\Gamma(z)$ is Euler gamma function.
\end{Tm}
The operator $Q$ is an analogue of the diagonal matrix. Thus, formula \eqref{9.2} reduces  operator $A$ to the simplest form, namely, to the  diagonal form. It is interesting that
the operator $B$ satisfies the following relation
\begin{equation}B={\mathcal{J}}^{i\alpha},\label{9.5}\end{equation}
where ${\mathcal{J}}$ is the integration operator
\begin{equation}{\mathcal{J}}f=\int_{0}^{x}f(t)dt,\quad 0{<}x{<}\ell .\label{9.6}\end{equation}
Later, we considered  a broad class of triangular Liv\v{s}ic models \eqref{4.26} and proved \cite{Sak1} that  these models  are (under certain  conditions)
 linearly equivalent to the diagonal operators.
\section{Christoffel--Darboux  type formulas}
Let us consider operators $S_{\zeta}$ of the form
\begin{align} & S_{\zeta}f=f(x)-\int_{\zeta}^{\infty}k(x,t)f(t)dt,\quad \zeta{\geq}0,\label{10.2}
\\ &
k(x,t)=\frac{A(x)A^{\prime}(t)-
A^{\prime}(x)A(t)}{\pi(x-t)},\label{10.3}\end{align}
which generalize operators \eqref{6.46}, \eqref{6.47}. It is easy to see that \eqref{10.3} yields
the following lemma.
\begin{La}\label{Lemma 10.1} Let $A(x)$ satisfy the equation
\begin{equation}A^{\prime\prime}(x)=v(x)A(x),\quad 0<x<\omega{\leq}\infty,
\label{10.4}\end{equation}
where $v(x)$ is a continuous function. Then, we have
\begin{equation}\frac{\partial}{{\partial}x}k(x,t)+\frac{\partial}{{\partial}t}k(x,t)=
A(t)\frac{v(t)-v(x)}{\pi(x-t)}A(x).\label{10.5}\end{equation}\end{La}

Introduce the notation
\begin{equation}k_{1}(\zeta,\eta)=k\left(\frac{\zeta+\eta}{2},\frac{\zeta-\eta}{2}\right).\label{10.6}
\end{equation}
After the change of  variables $\zeta=x+t$ and $\eta=x-t$, we rewrite \eqref{10.5} in the form
\begin{equation}2\frac{\partial}{\partial\zeta}k_{1}(\zeta,\eta)=
A\left(\frac{\zeta-\eta}{2}\right)\frac{v(\frac{\zeta-\eta}{2})-v(\frac{\zeta+\eta}{2})}{\pi\eta}A\left(\frac{\zeta+\eta}{2}\right).\label{10.7}\end{equation}
Taking into account relation \eqref{10.7}, we  obtain the next assertion.
\begin{La}\label{Lemma 10.2} Let the conditions of Lemma \ref{Lemma 10.1} be fulfilled and assume that
\begin{equation}\left|A\left(\frac{\zeta-\eta}{2}\right)\left[v\left(\frac{\zeta-\eta}{2}\right)-v\left(\frac{\zeta+\eta}{2}\right)\right]A\left(\frac{\zeta+\eta}{2}\right)\right|{\in}
L(\zeta,\infty),\label{10.8}\end{equation} where $0{\leq}\zeta<\infty$ and $-\zeta{\leq}\eta{\leq}\zeta.$ Then, we have
\begin{equation}
k_{1}(\zeta,\eta)=-\frac{1}{2}\int_{\zeta}^{\infty}
A\left(\frac{s-\eta}{2}\right)\frac{v(\frac{s-\eta}{2})-v(\frac{s+\eta}{2})}{\pi\eta}A\left(\frac{s+\eta}{2}\right)ds+g(\eta),\label{10.9}\end{equation}
where $g(\eta)$ is some differentiable function.\end{La}
Lemma \ref{Lemma 10.1} and Lemma \ref{Lemma 10.2} imply the assertion:
\begin{Pn}\label{Proposition 10.3} Let the conditions of Lemmas  \ref{Lemma 10.1} and \ref{Lemma 10.2} be fulfilled. If
\begin{equation}\lim_{\zeta{\to}+\infty}{k_{1}(\zeta,\eta)}=0,\label{10.10-}\end{equation}
then $g(\eta)$ in the right-hand side of \eqref{10.9} equals zero. \end{Pn}
\begin{Cy}\label{Corollary 10.4} Let the conditions of Proposition \ref{Proposition 10.3} be fulfilled.
Coming back  to the variables $x$ and $t$ in \eqref{10.9}, we have
\begin{equation}
k(x,t)=-\frac{1}{2}\int_{x+t}^{\infty}
A\left(\frac{s-x+t}{2}\right)\frac{v(\frac{s-x+t}{2})-v(\frac{s+x-t}{2})}{\pi(x-t)}A\left(\frac{s+x-t}{2}\right)ds\label{10.10}\end{equation}
or, equivalently,
\begin{equation} k(x,t)=\int_{0}^{\infty}
A(u+x)\frac{[v(u+x)-v(u+t)]}{\pi(x-t)}A(u+t)du.\label{10.12}\end{equation}\end{Cy}
\begin{Ee}\label{Example 10.4} Let us consider the case, where $v(x)=x.$\end{Ee}
In this case, we have $A(x)=\sqrt{\pi}A_{i}(x)$ where $A_{i}(x)$ is Airy function.  All the conditions of Proposition \ref{Proposition 10.3} are fulfilled (see \eqref{6.48}).
Therefore, using
\eqref{10.12} we obtain the well-known (see \cite{TrWi}) formula \eqref{6.50}.
\begin{Ee}\label{Example 10.5} Let us consider the case, where $v(x)=x^{2}.$ \end{Ee}
In this case, we have $A(x)=\sqrt{\pi}\exp(-\frac{1}{2}x^2).$ Using equality \eqref{6.47}, we obtain
\begin{equation}k(x,t)=\exp\left[-\frac{1}{2}(x^2+t^2)\right].\label{10.13}\end{equation}
Hence, the operator $S_{\zeta}$ has the form
\begin{equation}S_{\zeta}f=f(x)-\int_{\zeta}^{\infty}\exp\left[-\frac{1}{2}(x^2+t^2)\right]f(t)dt.
\label{10.14}\end{equation}
It follows from \eqref{10.14} that
\begin{equation}\det S_{\zeta}=1-\int_{\zeta}^{\infty}\exp(-t^2)dt.
\label{10.15}\end{equation}
Let us introduce the error  function erf:
\begin{equation}\erf(\zeta)=\frac{2}{\sqrt\pi}\int_{0}^{\zeta}\exp(-t^2)dt.\label{10.16}\end{equation}
Formulas \eqref{10.15} and \eqref{10.16} imply that
\begin{equation}\det S_{\zeta}=1-\frac{\sqrt\pi}{2}[1-\erf(\zeta)].\label{10.17}\end{equation} Using
the asymptotic expansion  of the error function for large real $\zeta$ (see \cite{Dw}),
we obtain
\begin{equation}\det S_{\zeta}=1-\frac{\exp(-\zeta^2)}{2x}\sum_{n=0}^{\infty}\left[(-1)^{n}
\frac{(2n-1)!!}{2^{n-1}\zeta^{2n}}\right].\label{10.18}\end{equation}
where (2n-1)!!  is the product of all positive odd numbers up to (2n-1) if n is positive and (-1)!!=1.

{\bf Acknowledgments.} The author is very grateful to A. Sakhnovich and I. Tydniouk for help and useful remarks.

\end{document}